# Quantitative indicators of the solutions of Diophantine equations and systems in the domain of the natural numbers

VICTOR VOLFSON

ABSTRACT. The paper shows that the asymptotic density of solutions of Diophantine equations or systems of the natural numbers is 0. The author provides estimation methods and estimates number, density and probability of $k$ - tuples $<x_1,...x_k>$ to be the solution of the algebraic equations of the first, second and higher orders in two or more variables, non-algebraic Diophantine equations and systems of Diophantine equations in the domain of the natural numbers. The estimate for the number of positive integer solutions of the second-order Diophantine equations in two, three or more variables is geometrically proved in the paper. The author proves the assertion about the number of solutions of algebraic Diophantine equations of higher orders in the domain of the natural numbers. The author provides the estimates for the asymptotic behavior of quantitative solutions of Diophantine equations and systems in the domain of the natural numbers.

1. INTRODUCTION

The central problem of Diophantine geometry is to determine the structure of the set of integer (positive) solutions of Diophantine equations and systems primarily to clarify the question whether this set is finite or not. In this respect, a basic heuristic assumption is that the system has usually a finite number of solutions if degree of the equations is much greater than number of variables [1].

General method of Diophantine approximation is well-known in the study of integer solutions of Diophantine equations [2]. However, this determines only a finite number of integer solutions of Diophantine equations in two variables and only for certain types of the corresponding plane curves. This method has been developed for the large number of variables using Schmidt's theorem about the simultaneous approximations [3].

__________________________

All above-mentioned results have the feature that they give only a qualitative picture of the set of integer solutions for the certain types of Diophantine equations and systems. However, they do not provide any quantitative estimates, which can be used to describe the set of these solutions. This article is focused on this issue.

Let us consider Diophantine equations and systems with regard to the number of their solutions in the domain of the natural numbers, density and asymptotic density of these solutions, the probability of a tuple of $k$ natural numbers to be the solution of Diophantine equations or systems in $k$ variables.

Here we define Diophantine equation in $k$ variables as the equation

$$F(x_1,...x_k) = 0, \tag{1.1}$$

where all variables take values in the domain of the natural numbers at the same time.

In this case, the system of Diophantine equations is an expression

$$F_1(x_1,...x_k) = 0;...F_m(x_1,...x_k) = 0, \tag{1.2}$$

where $F_i(x_1,...x_k) = 0$ is $i$-th Diophantine equation of the form (1.1) and $1 < m < k$.

The following assertions are proved in [4].

Assertion 1. The algebra $<\{A^k, \emptyset, P(A^k)\}, \cup, \setminus>$ is the algebra of the events where $P(A^k)$ is the set of all subsets of the non-empty set $A^k$, the set $A^k$ is the $k$-fold direct product of the set $A = 1,2,...N$, and $\cup, \setminus$ are subset operations of union and complement, respectively.

Assertion 2. The density of $k$ tuples on the set $B \subseteq A^k$ defined as

$$P(B) = \pi(B)/N^k \tag{1.3}$$

is the probability measure on the finite space of the events considered in assertion 1 and $\pi(B)$ is the number of solutions of Diophantine equations in the set of natural numbers $B$.

The asymptotic density of $k$-tuples on the set $B \subseteq A^k$ defined by the formula:

$$P'(B) = \lim_{N \to \infty} \frac{\pi(B)}{N^k} \tag{1.4}$$

is also considered in [4].

Denote by $B_N$ the set $B$ at a fixed value $N$ and denote by $\pi(B_N)$ the number of $k$-tuples $<x_1,...x_k>$ on the set $B_N$.

The sequence $\pi(B_N)$ is non-decreasing by definition.

Consider three cases for the asymptotic density of this sequence.

1) The sequence $\pi(B_N)$ does not increase starting from a certain value $N_0$, i.e. it remains constant. In this case, $P'(B) = \lim\limits_{N \to \infty} \dfrac{\pi(B_N)}{N^k} = 0$.

2) The sequence $\pi(B_N)$ increases without bound as $O(N^s)$ where $s < k$. In this case $P'(B) = \lim\limits_{N \to \infty} \dfrac{\pi(B_N)}{N^k} = 0$.

3) The sequence $\pi(B_N)$ increases without bound as $O(N^s)$ where $s = k$ but $\pi(B_N) \leq N^k$. In this case $P'(B) = \lim\limits_{N \to \infty} \dfrac{\pi(B_N)}{N^k} = a$, where the condition $0 < a \leq 1$ is satisfied for the constant $a$.

Now we consider Diophantine equations and systems using the above parameters and classification of the sequence $\pi(B_N)$ depending on the asymptotic density in the cases mentioned above.

2. DIOPHANTINE EQUATIONS WITH THE EXPLICIT VARIABLE

If any of the variables in Diophantine equation (1.1) is explicitly expressed, this equation can be written as

$$x_i = f(x_1,...x_{i-1}, x_{i+1},...x_k), \qquad (2.1)$$

where $f$ is the function taking values in the domain of the natural numbers. In (2.1) the variables $x_1,...x_{i-1}, x_{i+1},...x_k$ are independent and can take any value from the domain of natural numbers in the domain of the function $f$ and $x_i$ is dependent variable taking values of the function $f$ only.

Assume that $k-1$ independent variables in (2.1) take values in $A$, where $A = 1, 2, ... N$. Then, it follows from the definition of the function $f$ that the number of these values does not exceed

$$N^{k-1}. \tag{2.2}$$

Since the function $f$ is the surjective, the number of solutions of Diophantine equation (2.1) in the domain $A^k$ does not exceed $N^{k-1}$, i.e.

$$\pi(B_N) \leq N^{k-1}. \tag{2.3}$$

If the function $f$ is the injective, then the number of solutions of Diophantine equation (2.1) in the domain $A^k$ is $N^{k-1}$, as the function is one-to-one mapping in this case.

The equation $x_2 = x_1$ is an example of this Diophantine equation. The solutions of Diophantine equation lie on the main diagonal of the first quadrant: $(1,1); (2,2); ...(N,N)$, hence $A^2$ contains exactly $N$ solutions, which corresponds to the formula (2.2) for $k = 2$.

If the function $f$ is not injective in the domain $A^k$, then the number of solutions of the Diophantine equation (2.1) is less than $N^{k-1}$.

Based on (1.3), (2.3), the probability of natural $k$-tuples $<x_1,...x_k>$ to be solutions of Diophantine equation (2.1) in the domain $A^k$ (equal to the density of solutions of this equation) does not exceed:

$$Pr(B_N) = P(B_N) = \pi(B_N) / N^k \leq N^{k-1} / N^k = 1/N. \tag{2.4}$$

The following assertion results from the formula (2.4).

Assertion 3. The asymptotic density of the solutions of Diophantine equation (2.1) is equal to 0.

Thus, the asymptotic density of the solutions of Diophantine equation (2.1) refers to cases 1 and 2 for the asymptotic density discussed earlier. It is interesting that the asymptotic density of the solutions of Diophantine equation (2.1) is significantly lower than that of the coprimes, which refers to the case 3 [4].

Now consider Diophantine equation (1.1) in the case where the variables of the equation can't be explicitly expressed.

## 3. ALGEBRAIC DIOPHANTINE EQUATIONS

We begin our consideration of Diophantine equations with algebraic equations in one unknown. Consider the equation: $a_0 x^m + a_1 x^{m-1} + ... + a_{m-1} x + a_m = 0$, (3.1) where $m$ is a natural number and $a_0, ... a_m$ are integer and nonzero.

As all positive solutions of (3.1) are known to be divisors of the constant term, the number of such solutions is finite. Consequently, the asymptotic density of solutions of this equation is equal to 0 (the case 1 of the asymptotic density).

Let us consider the linear Diophantine equation of any number of variables. In [5] divisibility of $b$ by the greatest common divisor of the coefficients $a_1, a_2, ... a_k$ is shown to be necessary and sufficient condition for the linear equation

$$a_1 x_1 + a_2 x_2 + ... + a_k x_k = b \quad (3.2)$$

(where $a_1, a_2, ... a_k, b$ are integer) to be solvable in integers.

Now we consider the equation

$$a_2 x_2 - a_1 x_1 = b \quad (3.3)$$

as a special case of (3.2) where $a_1, a_2$ and $b$ are natural numbers. If equation (3.3) is solvable in integers, then the number $b$ is divisible by the greatest common divisor $a_1, a_2$. Then, if $x_{10}, x_{20}$ are the solutions of equation (3.3), we have that $a_2(x_{20} + a_1 t) - a_1(x_{10} + a_2 t) = b$ for any positive integer $t$. Since $a_1, a_2$ are natural numbers and sequences $x_{20} + a_1 t, x_{10} + a_2 t$ are strictly increasing, then the solutions of (3.3) are the natural numbers starting from some t. It is also given that $a_1, a_2 \geq 1$. Then the number of solutions of equation (3.3) is less than $N$ in the domain $A^2$ (where $A = 1, 2, ..N$) and the asymptotic density of the solutions is 0 (the case 2 of the asymptotic density).

In another way, we solve linear Diophantine with two variables:

$$a_1 x_1 + a_2 x_2 = b, \quad (3.4)$$

where $a_1, a_2$ and $b$ are natural numbers.

Assume that equation (3.4) has a solution in integers. It means that $b$ is divisible by the greatest common divisor $a_1, a_2$ - d. Divide equation (3.4) by d and obtain the equation with coprimes $a_1, a_2$. Continue to consider the equation (3.4) with coprime coefficients. If $b = a_1 \cdot a_2$, then equation (3.4) has no positive integer solutions, because otherwise $a_1 x_1 + a_2 x_2 = a_1 a_2$ and therefore $a_1 x_1 = a_2(a_1 - x_2)$. Since the numbers $a_1, a_2$ are coprimes, then $x_1$ is divisible by $a_2$ and we obtain $x_1 \geq a_2$; therefore, $a_1 x_1 + a_2 x_2 \geq a_1 a_2$, that contradicts to the fact that $a_1 x_1 + a_2 x_2 = a_1 a_2$.

Equation (3.4) has a finite number of solutions in positive integers if $b > a_1 a_2$, as in this case there are natural numbers $u, v$ such that the solution of equation (3.4) can be written as: $x_1 = u - a_1 t, x_2 = a_1 t - v$ where $t = 0, 1, 2.....$ (the case 1 of the asymptotic density).

In the general case, equation (3.2) has a finite number of natural solutions, if $a_1, ... a_k, b$ are natural numbers and $b$ is divided by the greatest common divisor of the coefficients $a_1, ... a_k$ as $x_i \leq b$ is satisfied (the case 1 of the asymptotic density).

Summarize the consideration of linear equation (3.2) $a_1 x_1 + ... + a_k x_k = b$ where the coefficients $a_i$ and $b$ are integers:

1. If the number $b$ is not divisible by the greatest common divisor $a_1, ... a_k$, then equation (3.2) has no integer solutions and therefore has no natural solutions (the case 1 of the asymptotic density).

2. If $a_1, ... a_k, b$ are positive integers (natural numbers) and the number $b$ is divided by the greatest common divisor of the coefficients $a_1, ... a_k$, then equation (3.2) has a finite number of positive solutions (the case 1 of the asymptotic density).

3. If any of the coefficients of the equation $a_1, ... a_k$ is negative integer, then the equation can have an infinite number of solutions (the case 2 of the asymptotic density). Equation (3.3) is an example of this case.

Now let us consider the nonlinear algebraic equations, and in particular Diophantine equation of the second order in two variables.

In [6] we give a solution of the integer-valued equation $Ax^2 + Bxy + Cy^2 = N$ in integers and find the exact number of this solution on the basis of the theory of comparisons. In this paper, we are only interested to estimate the number of solutions of the Diophantine equations in the domain of the natural numbers. We use the less accurate, but more visual geometric proofs to obtain these estimates.

We will proceed from the fact that the solutions of any Diophantine equation of the second order in two variables in integer (natural) numbers are integer (natural) coordinates on a curve of the second order:

$$F(x_1, x_2) = a_{11}x_1^2 + a_{22}x_2^2 + 2a_{12} + 2a_{13}x_1 + 2a_{23}x_2 + a_3 3 = 0 \tag{3.5}$$

where all $a_{ij}$ are integers.

Begin considering with the cases of non-degenerate curves of the second order.

There are three types of non-singular curves of the second order: ellipse (imaginary ellipse), hyperbole and parabola [7].

The determinant

$$D = a_{11}a_{22} - a_{12}^2 \tag{3.6}$$

is an invariant used to determine the type of non-degenerate second-order curve.

Assertion 4. The second-order Diophantine equation in two variables has a finite number of integer (positive) solutions where the equation $F(x_1, x_2) = 0$ defines any ellipse ($D = a_{11}a_{22} - a_{12}^2 > 0$) (the case 1 of the asymptotic density).

Proof. Plot the rectangular grid with integer nodes $x_1, x_2$ on a plane. The straight lines of the coordinate grid perpendicular to the axis $x_1$ and parallel to the axis $x_2$ cross an ellipse in a finite number of points with integer coordinate $x_1$. Naturally, $x_2$ is also integer for not all the points of intersection of the ellipse, so the numbers of points on the ellipse where the two coordinates are integer are finite.

Note. The imaginary ellipse did not have any real point. Therefore, the second order Diophantine equation in two variables has no solutions in integer (natural) numbers where the equation $F(x_1, x_2) = 0$ defines the imaginary ellipse.

Assertion 5. The second-order Diophantine equation in two variables can have an infinite number of solutions in the integer (natural) numbers where the equation $F(x_1, x_2) = 0$ defines a hyperbola or parabola.

Proof. Plot the rectangular grid with integer nodes $x_1, x_2$ on a plane. The straight lines of the coordinate grid perpendicular to the axis $x_1$ and parallel to the axis $x_2$ cross the hyperbola or parabola in an infinite number of points where the coordinate $x_1$ is integer. Therefore, if an infinite number of these points have integer coordinate $x_2$, then equation (3.5) will have an infinite number of solutions in the integer (natural) numbers.

We remind that the invariant

$$D = a_{11}a_{22} - a_{12}^2 < 0 \qquad (3.7)$$

for any hyperbole.

The invariant

$$D = a_{11}a_{22} - a_{12}^2 = 0 \qquad (3.8)$$

for any parabola.

Therefore the assertion 7 can be formulated on the basis (3.7) and (3.8) in another way. The second-order Diophantine equation in two variables (3.5) can have an infinite number of the integer (natural) solutions if the condition

$$D = a_{11}a_{22} - a_{12}^2 \leq 0 \qquad (3.9)$$

is held for its coefficients.

It will be shown in assertion 14 that the number of solutions of algebraic equation (1.1) generally does not exceed $nN^{k-1}$ in $A^k$ (the case 2 of the asymptotic density).

Now we consider the second order Diophantine equation in two variables where the equation (3.5) defines the second-order degenerate curves. Let us start with a degenerate ellipse. The degenerate ellipse is a real point, so if there exists a solution in integer (natural) numbers, it should be unique one. Thus, in this case, as well as for non-degenerate ellipse, the number of solutions is finite (the case 1 of the asymptotic density).

We remind that the invariant $D > 0$ for non-degenerate ellipse, as for an ellipse, and it is enough to check this condition.

One can easily show that in the case of a degenerate ellipse the solution is expressed in real numbers in terms of the coefficients of equation (3.5) as follows:

$$(-a_{13}/a_{11}, -a_{23}/a_{22}). \tag{3.10}$$

Hence, it is clear what is required to ensure the integer (natural) solution.

Now we consider the case when equation (3.5) $F(x_1, x_2) = 0$ defines a degenerate hyperbola. A degenerate hyperbola consists of the two intersecting lines symmetric with respect to the point of their intersection. One line is increasing and another one is decreasing function of $x_1$.

There can exist only a finite number of solutions for equation (3.5) on the decreasing line and an infinite number of solutions on the increasing line in the domain of the natural numbers. The number of solutions does not exceed $N$ on the increasing line in the domain $A^2$ where $A = 1, 2...N$.

In the case of a degenerate hyperbola, the total number of the solutions of equation (3.5) does not exceed $2N$ on both lines (the case 2 of the asymptotic density). The number of the solutions is equal to $2N$ in the domain $A^2$ if the lines are parallel to the coordinate axes (the case 2 of the asymptotic density).

Consider the case when equation (3.5) $F(x_1, x_2) = 0$ defines a degenerate parabola. This case is divided into three subcases.

1. Two parallel lines. Straight lines are parallel to either the axis $x_1$ or the axis $x_2$, so the solutions are of the form $x_1 = b_1, x_1 = b_2$ or $x_2 = c_1, x_2 = c_2$. If numbers $b_1, b_2$ or $c_1, c_2$ are natural, then equation (3.5) has an infinite number of natural solutions in the domain c (the case 2 of the asymptotic density).

2. One straight line is parallel to either the axis $x_1$ or the axis $x_2$, so the solution has the form $x_1 = a$ or $x_2 = b$. If $a, b$ are the natural numbers, then equation (3.5) has an infinite number of solutions in the domain of the natural numbers, and $N$ solutions in the domain $A^2$ (the case 2 of the asymptotic density).

3. The imaginary parallel lines. In this case, equation (3.5) has no solutions in the real numbers and, respectively, in the integer (natural) numbers (the case 1 of the asymptotic density).

The second-order Diophantine equation (1.1) in three variables can be represented as:

$$F(x_1,x_2,x_3) = a_{11}x_1^2 + a_{22}x_2^2 + a_{33}x_3^2 + a_{12}x_1x_2 + a_{13}x_1x_3 + a_{23}x_2x_3 + 2a_{14}x_1 + 2a_{24}x_2 + 2a_{34}x_3 + a_{44} = 0 \quad (3.11)$$

Equation (3.11) corresponds to the second-order surface. There are 17 types of the second-order surfaces [7]. I recall that the second-order surfaces have their invariants, which do not change under the orthogonal transformation:

$$I = a_{11} + a_{22} + a_{33}, \tag{3.12}$$

$$A' = A_{11} + A_{22} + A_{33} + A_{44}, \tag{3.13}$$

where $A_{ij}$ is the algebraic addition to the $a_{ij}$ in $D$.

$$|J| = \begin{bmatrix} a_{11} & a_{12} \\ a_{12} & a_{22} \end{bmatrix} + \begin{bmatrix} a_{22} & a_{23} \\ a_{23} & a_{33} \end{bmatrix} + \begin{bmatrix} a_{33} & a_{13} \\ a_{13} & a_{11} \end{bmatrix} \tag{3.14}$$

$$|D| = \begin{bmatrix} a_{11} & a_{12} & a_{13} \\ a_{12} & a_{22} & a_{23} \\ a_{13} & a_{23} & a_{33} \end{bmatrix}, \tag{3.15}$$

$$|A| = \begin{bmatrix} a_{11} & a_{12} & a_{13} & a_{14} \\ a_{12} & a_{22} & a_{23} & a_{24} \\ a_{13} & a_{23} & a_{33} & a_{34} \\ a_{14} & a_{24} & a_{34} & a_{44} \end{bmatrix}. \tag{3.16}$$

We begin considering of the second-order equation in three variables with non-degenerate second-order surfaces, for which $(3.16) |A| \neq 0$.

Assertion 8. The second order Diophantine equation $F(x_1, x_2, x_3) = 0$ has the finite number of solutions where it defines an ellipsoid (the case 1 of the asymptotic density).

Proof. We project the ellipsoid into the plane $x_1, 0, x_2$. There is a finite number of points with integer coordinates within the projection (ellipse). Let us draw perpendiculars to the plane from the points with the integer coordinates to the intersection with the ellipsoid.

Not all the points of intersection will have integer coordinates, so the number of the points with integer (natural) coordinates will be finite.

The surface (3.11) is ellipsoid, if the conditions $(|A|<0, |D|\neq 0, I|D|>0)$ are satisfied as it follows from (3.12) - (3.16).

Note 1: Second-order Diophantine equation (3.11) $F(x_1, x_2, x_3) = 0$ has no solutions in integers (positive integers) where the equation defines an imaginary ellipsoid $(|A|>0, |D|\neq 0, I|D|<0)$ as the imaginary ellipsoid does not have any real point (the case 1 of the asymptotic density).

Note 2: These statements are held for the second-order equation in $k$ variables where the equation defines an ellipsoid (imaginary ellipsoid) in k-dimensional space of the rectangular coordinates (the case 1 of the asymptotic density).

Assertion 9. The second-order Diophantine equation in three variables (3.11) can have an infinite number of solutions in the domain of the integer (natural) numbers where the equation $F(x_1, x_2, x_3) = 0$ defines the non-degenerate quadric surface: the two-sheeted hyperboloid $(|A|<0, |D|\neq 0, I|D|<0)$, the one-sheeted hyperboloid $(|A|>0, |D|\neq 0, I|D|<0)$, the elliptic paraboloid $(|A|<0, |D|=0)$ and the hyperbolic paraboloid.

Proof. Assume without loss of generality that the cutting plane parallel to $x_1, 0, x_2$ contains either hyperbole or parabola. Project the surface onto the plane $x_1, 0, x_2$. There is an infinite number of points with integer (natural) coordinates $(x_1, x_2)$ within the projection of the hyperbola or parabola. Draw the perpendicular to the plane $x_1, 0, x_2$ from the points with integer coordinates $(x_1, x_2)$. Assume that these perpendiculars intersect the surface in an infinite number of points with the integer (natural) coordinate $x_3$. In this case, the number of integer (natural) solutions of equation (3.11) is infinite (the case 2 of the asymptotic density as it follows from Assertion 16).

Now consider the second-order Diophantine equations in three variables (3.11) where the equation $F(x_1, x_2, x_3) = 0$ defines degenerate quadric surface, for which the invariant $|A|=0$.

First, we consider the equations that define the conical surfaces: real $(|D|\neq 0, I|D|\leq 0)$ and imaginary cone $(|D|\neq 0, I|D|>0)$.

Assertion 10. The second-order Diophantine equation in three variables (3.11) can have an infinite number of solutions in the integer (natural) numbers where the equation $F(x_1, x_2, x_3) = 0$ defines a real cone.

Proof. Assume, without loss of generality, that the cutting planes parallel to $x_1, 0, x_2$ contain the intersecting lines. Project the surface to the plane $x_1, 0, x_2$. There exist an infinite number of points with integer (natural) coordinates within the intersecting lines $(x_1, x_2)$. Draw perpendiculars to the plane $x_1, 0, x_2$ from the points with integer coordinates $(x_1, x_2)$ and assume that the perpendiculars intersect the surface in an infinite number of points with integer (natural) coordinate $x_3$. In this case, the number of integer (natural) solutions of equation (3.11) is infinite (the case 2 of the asymptotic density as it follows from Assertion 16).

Note 3: This statement is true for the second-order equations in $k$ variables where the equation $F(x_1, x_2, ... x_k) = 0$ defines a real cone in $k$-dimensional space of rectangular coordinates.

The imaginary cone has unique real point, and therefore the equation can only have not more than unique integer (natural) solution in this case (the case 1 of the asymptotic density).

Now we consider the Diophantine equation of the form (3.11) where the equation $F(x_1, x_2, x_3) = 0$ is the cylindrical surface of the second order. There are six cylindrical surfaces of the second order: elliptic cylinder, parabolic cylinder, hyperbolic cylinder, a pair of coincident lines, a pair of coincident planes and a pair of intersecting planes [7]. Estimate the number of solutions of equations (3.11) for all types of the cylindrical surfaces of the second order.

Assertion 11. A number of the solutions of Diophantine equation (3.11) where it corresponds to the elliptic cylinder in the domain of the natural numbers $A^3, A = 1, 2...N$ is equal to $\pi(B_N) = kN$, the probability of a tuple $< x_1, x_2, x_3 >$ to be the solution of this equation is equal to

$$Pr(B_N) = k / N^2, \qquad (3.17)$$

where $k$ is the finite non-negative number (the case 2 of the asymptotic density).

Proof. Assume, without sacrificing generality, that the intersection of the elliptic cylinder with the plane $x_1, 0, x_2$ is an ellipse. As it follows from the assertion proved above, the equation where $F(x_1, x_2) = 0$ is the ellipse has only a finite number of solutions in positive integers ($k$). Draw perpendiculars to the planes $x_1, 0, x_2$ passing through the points of solutions of the equation $F(x_1, x_2) = 0$. There are $N$ solutions of Equation (33) on each perpendicular in the domain $A^3$ and $kN$ solutions on the $k$ perpendiculars respectively. Consequently, the probability of the tuple $<x_1, x_2, x_3>$ to be the solution of this equation is equal to - $Pr(B_N) = kN / N^3 = k / N^2$.

Assertion 12. The number of solutions of Diophantine equation (3.11) does not exceed $\pi(B_N) \leq N^2$ where Diophantine equation defines a parabolic or a hyperbolic cylinder in the domain $A^3$. Consequently, the probability of the tuple $<x_1, x_2, x_3>$ to be the solution of this equation does not exceed

$$Pr(B_N) \leq 1 / N \qquad (3.18)$$

(the case 2 of the asymptotic density).

Proof. Assume, without sacrificing generality, that the intersection of the parabolic or the hyperbolic cylinder with the plane $x_1, 0, x_2$ is a parabola or hyperbola. The equation $F(x_1, x_2) = 0$ can have an infinite number of solutions in the domain of the natural numbers if this equation is parabola or hyperbola on the basis of the earlier assertion proved. The number of solutions of this equation in the domain $A^2$ does not exceed $N$. Draw perpendiculars to the plane $x_1, 0, x_2$ passing through the points of solutions of the equation $F(x_1, x_2) = 0$. There are $N$ solutions of the equation (3.11) on each perpendicular in the domain $A^3$ and not more than $N^2$ solutions on $N$ perpendiculars, respectively. Consequently, $\pi(B_N) \leq N^2$ and the probability of the tuple $<x_1, x_2, x_3>$ to be the solution of this equation is equal to:

$$Pr(B_N) \leq N^2 / N^3 = 1 / N. \qquad (3.19)$$

Equation (3.11) has no solutions in real (and even more natural) numbers, if (3.11) defines an imaginary elliptical cylinder.

The solutions of (3.11) lie in a plane parallel to a coordinate plane $x_2, 0, x_3$, if (3.11) defines a pair of coincident planes. The number of solutions of this equation is equal to $N^2$ in the domain $A^3$ and the probability of a tuple $<x_1, x_2, x_3>$ to be the solution of this equation is equal to $Pr(B_N) = N^2 / N^3 = 1/N$.

Assertion 13. If (3.11) corresponds to a pair of intersecting planes, then the number of solutions of the equation $F(x_1, x_2, x_3) = 0$ does not exceed $\pi(B_N) = N^2 + kN$ in the domain $A^3$ where $k$ is non-negative integer. The probability of a tuple $<x_1, x_2, x_3>$ to be the solution of this equation does not exceed:

$$Pr(B_N) \leq 1/N + k/N^2 \qquad (3.20)$$

(the case 2 of the asymptotic density).

Proof. Assume, without loss of generality, that these planes are projected to the coordinate plane $x_1, 0, x_2$ as two intersecting lines. It was previously shown that the equation $F(x_1, x_2) = 0$ corresponded to the intersecting straight lines. Therefore, the number of its solutions does not exceed $N + k$ in the domain $A^2$. Draw the perpendicular lines to the plane $x_1, 0, x_2$ passing through the solutions of the equation $F(x_1, x_2) = 0$. Each perpendicular line contains $N$ solutions, so the number of solutions of equation (3.11) does not exceed $\pi(B_N) \leq N^2 + kN$ in the domain $A^3$ and the probability of a tuple $<x_1, x_2, x_3>$ to be the solution of this equation does not exceed:

$$Pr(B_N) \leq (N^2 + kN) / N^3 = 1/N + k/N^2 \qquad (3.21)$$

Assertion 14. The number of solutions of equation (3.11) does not exceed $\pi(B_N) \leq 2N^2$ in the domain $A^3$ if the equation corresponds to a pair of parallel planes. The probability of the tuple $<x_1, x_2, x_3>$ to be the solution of this equation $F(x_1, x_2, x_3) = 0$ does not exceed

$$Pr(B_N) \leq 2/N \qquad (3.22)$$

(the case 2 of the asymptotic density).

Proof. Assume, without loss of generality, that these planes are projected onto the coordinate plane $x_1, 0, x_2$ as two parallel lines. It was previously shown that the number of solutions of equation (3.11) did not exceed $2N$ in the domain $A^2$ if $F(x_1, x_2) = 0$ corresponds to the parallel straight lines. Draw perpendicular lines to the plane $x_1, 0, x_2$ passing through the solutions of the equation $F(x_1, x_2) = 0$.

Each perpendicular line contains $N$ solutions, so the number of solutions (3.11) does not exceed $\pi(B_N) \leq 2N^2$ in the domain $A^3$, and the probability of a tuple $<x_1, x_2, x_3>$ to be the solution of this equation does not exceed $Pr(B_N) \leq (2N^2)/N^3 = 2/N$ that corresponds to (3.22).

Finally, if Diophantine equation (3.11) corresponds to the pair of imaginary parallel planes, then the equation $F(x_1, x_2, x_3) = 0$ has no solutions in real numbers and all the more in the natural numbers.

We should take into account the domain of the function $F(x_1, ... x_k) = 0$, when considering the number of solutions of the Diophantine equation. The number of solutions in the domain of the natural numbers $A^k$ ($A = 1, 2, ... N$) may be reduced, if we consider the function domain.

Assertion 15. The number of solutions (1.1) in the domain $A^k$ does not exceed $N^{k-1}$, if there exists even one surjective map $x_i = f_i(x_1, ... x_{i-1}, x_{i+1}, ... x_k), 1 \leq i \leq k$ corresponding to equation (1.1) (the case 2 of the asymptotic density).

Proof. If there exists even one surjective map $x_i = f_i(x_1, ... x_{i-1}, x_{i+1}, ... x_k)$, then a unique value $x_i$ corresponds to each value of the variables $x_1, ... x_{i-1}, x_{i+1}, ... x_k$ for this map in the domain $A^k$.

If all the values of variables $x_1, ... x_{i-1}, x_{i+1}, ... x_k$ are included in the domain of function (1.1), then each of $k - 1$ variables takes $N$ values in the domain $A^k$. Since the values $x_i$ can be repeated, all variables take not more than $N^{k-1}$ values.

If any of the values of variables $x_1, ... x_{i-1}, x_{i+1}, ... x_k$ is not included in the domain of function (1.1), then some variable takes less than $N$ values in the domain $A^k$, i.e. all variables take less than $N^{k-1}$ values.

This assertion generalizes formula (7) not only for algebraic Diophantine equations, but also for the case where even one variable in equation (1.1) can be explicitly expressed in terms of the other variables in the domain $A^k$.

If all the maps $x_i = f_i(x_1,...x_{i-1}, x_{i+1},...x_k), 1 \leq i \leq k$ corresponding to the algebraic Diophantine equations of $n$-order (1.1) $F(x_1,...x_k) = 0$ are not surjective, then the number of solutions of (1.1) can exceed $N^{k-1}$ in the domain $A^k$ where $A = 1, 2,...N$.

Assertion 16. The number and the density of solutions of the algebraic Diophantine equation of n-order (1) does not exceed $\pi(B_N) \leq nN^{k-1}$ and $P(B_N) \leq n/N$, respectively, in the domain $A^k$ where $A = 1, 2,...N$ and the asymptotic density is equal to 0 (the case 2 of the asymptotic density).

Proof. In the assertion 15 we considered the case where there exists even one surjective map. Now we consider the case where each map $x_i = f_i(x_1,...x_{i-1}, x_{i+1},...x_k), 1 \leq i \leq k$, corresponding to the $n$-order algebraic Diophantine equation $F(x_1,...x_k) = 0$ is not surjective.

The straight line may cross the $n$-order algebraic surface in $k$-dimensional space at no more than $n$ points, if it is not a generatrix. Each perpendicular drawn to the coordinate plane from the points with natural values of the coordinates intersects the surface at no more than $n$ points, if it is not a generatrix. There are $N^{k-1}$ perpendiculars in the domain $A^k$, therefore the maximum number of such intersections will be $nN^{k-1}$. Therefore, the number of solutions of the $n$-order algebraic Diophantine equation $F(x_1,...x_k) = 0$ does not exceed $nN^{k-1}$ in the domain $A^k$ where $A = 1, 2,...N$.

It can be shown that if rectilinear generatrices parallel to each coordinate axis pass through any point of the surface in $k$-dimensional space with Cartesian coordinates corresponding to the equation $F(x_1,...x_k) = 0$, then this surface is decomposed into the $k-1$ hyperplanes parallel to the coordinate axes. In this case, the number of solutions of this equation in positive integers does not exceed $(k-1)N^{k-1}$ in the domain $A^k$. From the relation $n > k-1$ it follows that $(k-1)N^{k-1} < nN^{k-1}$. Otherwise, a perpendicular passing through a point of the surface is not a generatrix in even one direction what was previously discussed.

The density of solutions of the $n$-order algebraic Diophantine equation $F(x_1,...x_k)=0$ does not exceed $P(B_N) \leq nN^{k-1}/N^k = n/N$ in the domain $A^k$ where $A = 1,2,...N$.

The asymptotic density of the solutions of the n-order algebraic Diophantine equation $F(x_1,...x_k)=0$ is $P'(B) = \lim_{N\to\infty} n/N = 0$ in the domain $A^k$ where $A = 1,2,...N$.

Assertion 17. Let there be given the Diophantine equation

$$F(x_1,...x_k)=0, \qquad (3.23)$$

where $F(x_1,...x_k)$ is the polynomial with integer coefficients of the unknowns and the free term, then:

1. Equation (3.23) has no solutions in positive integers, if all the coefficients of the polynomial are positive or equal to zero except for one and the constant term is not negative.

2. Equation (3.23) has a finite number of solutions in positive integers, if this equation corresponds to a non-degenerate surface and all the coefficients of the polynomial are positive or equal to zero except for one and the constant term is negative. The case the equation has no solutions (the number of solutions is 0) is related to a finite number of solutions.

3. Equation (3.23) can have an infinite number of solutions in positive integers, if the coefficients of the unknowns have different signs.

Proof.

1. a) The constant term is positive. $F(0,...0)$ is equal to the constant term and therefore $F(0,...0) > 0$. Since the value $F(x_1,...x_k)$ increases as a function of variables $x_1,...x_k$, then equation (3.23) has no solutions in the natural numbers of the variables.

b) Equation (3.23) has the solution $x_1 = 0,...x_k = 0$, if the free term is equal to 0 and $F(0,...0) = 0$. Since the value $F(x_1,...x_k)$ increases as the function of variables $x_1,...x_k$, then equation (3.23) has no solutions in the domain of the natural numbers.

2. $F(0,...0) < 0$ as the constant term is negative. The function $F(x_1,...x_k)$ increases by any variable and is equal to 0 on the surface with a finite value of the coordinates. Therefore, this surface may contain only a finite number of the solutions with natural coordinates and the corresponding equation (3.23) has a finite number of the solutions in the domain of the natural numbers.

3. It is enough to give an example to prove the assertion in this case.

There be given the equation:

$$x_2^m - x_1^n = 0, \qquad (3.24)$$

where $m, n$ and $m > n > 1$ are natural numbers. It is necessary to determine the number of solutions in the domain $A^2$.

Solution. Equation (3.24) has an infinite number of solutions in positive integers $x_1 = t^{l/n}, x_2 = t^{l/m}$ where $l$ is the least common multiple of positive integers $n, m$.

From $m > n$ it follows that $x_1 > x_2$ for $t \geq 2$. Thus, the number of solutions of (3.24) in the domain of the natural numbers $A^2$ is defined by the inequality $t^{l/n} < N$. Therefore

$$\pi(B_N) < [N^{n/l}], \qquad (3.25)$$

where $[\ ]$ is an integer part of the number.

Now consider the algebraic Diophantine equations that correspond to degenerate curves and surfaces of higher orders.

Let us ask ourselves, why does the equation $F(x_1, x_2) = 0$ corresponding to degenerate curve have more solutions than the equation corresponding to non-degenerate curve in the domain of the natural numbers $A^2$?

The answer is simple. The degenerate curve splits into the several lines. Each line can contain up to $N$ solutions in the domain $A^2$.

Maximum number of solutions is equal to

$$n \cdot N \qquad (3.26)$$

in the domain $A^2$ if the lines are parallel and their number is equal to $n$. The number of parallel lines, into which the curve splits, depends on its order. The $n$-order degenerate curve can contain $n$ parallel lines according to the greatest possible number of real roots for an algebraic equation of degree $n$. In this case, maximum number of solutions of the equation $F(x_1, x_2) = 0$ in the domain of the natural numbers $A^2$ is determined by formula (3.26).

If $F(x_1,...x_k) = 0$ corresponds to a degenerate surface in $k$-dimensional space, the solutions of Diophantine equations may be contained in several intersecting or parallel hyperplanes.

There can be no more than $N^{k-1}$ solutions in the domain $A^k$ of natural numbers. If the number of parallel hyperplanes is equal to $n$, maximum number of solutions of the equation $F(x_1,...x_k) = 0$ in the domain $A^k$ is equal to

$$n \cdot N^{k-1} \tag{3.27}$$

The $n$ order degenerate surface can contain $n$ parallel hyperplanes according to the greatest possible number of real roots of $n$ order algebraic equation. In this case, maximum number of solutions of the equation $F(x_1,...x_k) = 0$ in the domain of the natural numbers $A^k$ is defined by formula (3.27).

## 4. SYSTEMS OF DIOPHANTINE EQUATIONS

The system of Diophantine equations $F_1(x_1,...x_k) = 0,...F_m(x_1,...x_k) = 0$ (1.2) where $m < k$ is a special case of the predicate condition where the equations $F_i(x_1,..x_k) = 0$ (1.1) ($1 < i \leq m$) are associated with logical operation "AND". Thus, the solution of system (1.2) is understood as the intersection of the sets of solutions of equations (1.1).

Recall that the number of solutions of (1.1) is not greater than $nN^{k-1}$ in the domain $A^k$ for algebraic surfaces. Therefore, system (1.2) of Diophantine equations has no more than

$$\pi(B_N) \leq nN^{k-1} \tag{4.1}$$

solutions

As follows from (4.1), the probability of a tuple $<x_1,...x_k>$ to be the solution of the system of Diophantine equations (1.2) does not exceed

$$Pr(B_N) \leq n/N. \tag{4.2}$$

It follows from (4.2) that the asymptotic density of the solutions of system (1.2) is equal to

$$P'(B_N) = 0. \tag{4.3}$$

Based on (4.3), we can generalize Assertion 3 to the system of Diophantine equations, i.e. the asymptotic density of solutions of Diophantine equations and systems of Diophantine equations is equal to $0$.

The number of solutions of Diophantine equations is equal to the number of the solutions of an equation if only all the equations of system (1.2) are equivalent in the domain of the natural numbers.

The system of Diophantine equations (1.2) have no solutions in the natural numbers, if even one equation included in the system has no solutions in the natural numbers or the equations of the system have no common solutions in the natural numbers.

As an example, we find a number and density of solutions for the system of Diophantine equations:

$$x_1^{n_1} = x_2^{n_2} = \ldots = x_k^{n_k} \tag{4.4}$$

in the domain $A^k$ where $n_1, \ldots n_k$ are natural numbers and $n_k > n_{k-1} > \ldots > n_1$. System of equations (4.4) has infinitely many solutions in positive integers:

$$x_1 = t^{m/n_1}, x_2 = t^{m/n_2}, \ldots x_k = t^{m/n_k}, \tag{4.5}$$

where $m$ is the least common multiple of the numbers $n_1, n_2, \ldots n_k$, and $t$ is natural number. As $n_k > n_{k-1} > \ldots > n_1$, the number of solutions of the equation in the domain $A^k$ is determined by the inequality $t^{m/n_1} \leq N$ for $t > 2$. Therefore, the number of solutions of (4.5) is equal to $\pi(B_N) = [N^{n_1/m}] < N$ in the domain $A^k$. The density of solutions (4.4) is equal to $\pi(B_n) = [N^{n_1/m}]/N^k < 1/N^{k-1}$ in the domain $A^k$ (the case 2 of the asymptotic density).

5. NON-ALGEBRAIC DIOPHANTINE EQUATIONS

The simplest non-algebraic Diophantine equations are

$$x_2 = a^{x_1}, \tag{5.1}$$

$$x_2 = \log_a(x_1), \tag{5.2}$$

where $a$ is a natural number greater than 1.

Equations (5.1) and (5.2) have an infinite number of solutions in the natural numbers and $\log_a(N) < N$ solutions in the domain $A^2$.

It is easy to understand that the simplest Diophantine trigonometric and hyperbolic equations have no solutions in the natural numbers. The non-algebraic Diophantine equation can be represented as a system of simple non-algebraic equations and algebraic equations by introducing dummy variables.

For example, the equation

$$2^{x_1} + 3^{x_2} = 5^{x_3} \tag{5.3}$$

can be written as the system of equations

$$x_4 = 2^{x_1}, x_5 = 3^{x_2}, x_6 = 5^{x_3}, x_6 = x_4 + x_5. \tag{5.4}$$

The last equation of system (5.4) is algebraic. It has an infinite number of solutions and less than $N^2$ solutions in the domain $A^3$ (the case 2 of the asymptotic density).

The solutions of the systems of equations are the natural values of the variables that satisfy all the equations simultaneously and it is naturally that the system has fewer solutions than its any algebraic equation.

Therefore, Diophantine equations composed of trigonometric and hyperbolic functions have no solutions in the natural numbers.

Thus, Assertion 3 can be extended to the non-algebraic Diophantine equations.

For the system of the equations (5.3) the number of solutions is finite and equal to 2 as for equation (5.4) (the case 1 of the asymptotic density).

There are other simple Diophantine non-algebraic equations, for example,

$$x_3 = x_1^{x_2}, \tag{5.5}$$

where $x_1, x_2, x_3$ are the natural numbers and $x_1 > 1$.

Equation (5.5) has an infinite number of solutions in positive integers and $\sum_{i=2}^{N} \log_i(N) < N^2$ solutions in the domain $A^2$ (the case 2 of the asymptotic density).

Consider the equation

$$x_1^{x_2} - 2^{x_3} = 1. \tag{5.6}$$

Equation (5.6) is transformed into the system of equations consisting of two simple non-algebraic equations and the algebraic equations:

$$x_4 = x_1^{x_2}, x_5 = 2^{x_3}, x_4 - x_5 = 1. \tag{5.7}$$

The last algebraic equation in (5.7) has an infinite number of solutions in positive integers and less than $N$ solutions in the domain $A^2$ (the case 2 of the asymptotic density). However, to satisfy additional elementary equations, system (5.7) and equation (5.6) have only one solution in the domain of the natural numbers: $x_1 = 3, x_2 = 2, x_3 = 3$ (the case 1 of the asymptotic density).

6. ASYMPTOTIC BEHAVIOR OF THE QUANTITATIVE INDICATORS FOR THE SOLUTIONS OF DIOPHANTINE EQUATIONS AND SYSTEMS

Number, density and probability of $k$-tuples $<x_1,...x_k>$ to satisfy certain conditions in the domain $A^k$ (where $A = 1,...N$) determines the asymptotic behavior of these indicators, i.e. the functions $f(N), g(N)$, for which the corresponding asymptotic equality is satisfied:

$$\pi(B_N) \sim f(N), \tag{6.1}$$

$$P(B_N) = Pr(B_N) \sim g(N). \tag{6.2}$$

Example. Select a random tuple of the natural numbers $<x_1, x_2, x_3>$. Raise the first number to the square, the second one to the cube, and the third one to the fourth power.

It is required to determine such probability that all the three numbers are equal or smaller than $10^6$, i.e. $Pr(x_1^2 = x_2^3 = x_3^4 \leq 10^6)$. It is also required to determine quantity and density of the asymptotic behavior for tuples $<x_1, x_2, x_3>$ that satisfy the specified condition.

Solution. From (4.5) it follows that $\pi(B_N) = [N^{1/6}] = 10$ and therefore from (6.1) it follows that $\pi(B_N) \sim N^{1/6}$.

From (6.2) we have $P(B_N) = Pr(B_N) = [N^{1/6}]/N^3 = 10/10^{18} = 10^{-17}$, therefore, $P(B_N) = Pr(B_N) \sim N^{1/6}/N^3 = N^{-17/6}$.

Here is another example of how to determine the asymptotic behavior of quantitative solutions of Diophantine equation.

It is required to determine quantity and density of the asymptotic behavior of the solutions of Diophantine equations

$$2x_1^2 + x_2^2 - x_3^2 = 0 \qquad (6.3)$$

in the domain $A^2$. It is also required to determine the probability of a tuple $<x_1, x_2>$ to be the solution of Diophantine equation (6.3).

Since the number of solutions of (6.3) in the domain $A^2$ is equal to $\pi(B_N) = [\sqrt{(N-1)/2}]$, then

$$\pi(B_N) \sim \sqrt{N/2}. \qquad (6.4)$$

From (6.4) it follows that the asymptotic density of the solutions of (6.4) in the domain $A^2$ and the probability of a tuple $<x_1, x_2>$ to be the solution of Diophantine equation (6.4) is equal to: $P(B_N) = Pr(B_N) \sim \sqrt{N/2}/N^3 = 1/\sqrt{2N^5}$.

## 7. CONCLUSIONS AND SUGGESTIONS FOR FURTHER WORK

It is shown that the asymptotic density of solutions of Diophantine equations or systems in the natural numbers is 0.

The paper shows that the asymptotic density of solutions of Diophantine equations or systems of the natural numbers is 0. The author provides estimation methods and estimates number, density and probability of $k$ - tuples $<x_1,...x_k>$ to be the solution of the algebraic equations of the first, second and higher orders in two or more variables, non-algebraic Diophantine equations and systems of Diophantine equations in the domain of the natural numbers.

The author is going to consider the asymptotic density of the solutions of Diophantine inequities and systems in the natural numbers in the next paper.

## 8. ACKNOWLEDGEMENTS
Thank you, everyone, who has contributed to the discussion of this paper.

REFERENCES

1. Bashmakova I. G. Diofant and Diofant/s equations, Moscow, 1972. 68p.

2. S. Lang. Diophantine Geometry (3rd ed.) Interscience, New York (1962).

3. Borevich Z. I., Shafarevich I.R. Number Theory, Moscow, 1972.

4. Volfson V.L. Probabilistic properties of some sequences on a finite interval of the natural numbers Applied fhysics and mathematics, 2014. №6.

5. W. Sierpinski About solving equations in integers, Translated from the Polish IG Melnikova, M, Fizmatgiz, 1961. 88.

6. Buchstab A.A, Number Theory, Publishing House "Education", Moscow, 1966. 384 p.

7. Ilyin V.A, Pozdniak E.G. Analytical Geometry - M, Science,   1999. 224 p.